\documentclass[12pt]{amsart}
\usepackage[english]{babel}
\usepackage[T1]{fontenc}
\usepackage[utf8]{inputenc}
\usepackage{indentfirst}
\usepackage{fancyhdr}
\usepackage{amsthm}
\usepackage{amsmath}
\usepackage{amssymb}
\usepackage{braket}
\usepackage{enumitem}
\usepackage{color}
\usepackage{cite,enumitem,graphicx}
\usepackage[colorinlistoftodos]{todonotes}
\usepackage[colorlinks=true,urlcolor=blue,
citecolor=red,linkcolor=blue,linktocpage,
bookmarksnumbered,bookmarksopen,pdfpagelabels,]{hyperref}
\usepackage[hyperpageref]{backref}
\allowdisplaybreaks

\newtheorem{Th}{Theorem}[section]
\newtheorem{Prop}[Th]{Proposition}
\newtheorem{Lem}[Th]{Lemma}
\newtheorem{Cor}[Th]{Corollary}
\theoremstyle{definition}

\theoremstyle{remark}
\newtheorem{Rem}[Th]{Remark}

\newcommand{\R}{\mathbb{R}}
\newcommand{\cS}{\mathcal{S}}
\newcommand{\cD}{\mathcal{D}}
\newcommand{\cC}{\mathcal{C}}
\newcommand{\cM}{\mathcal{M}}
\newcommand{\cO}{\mathcal{O}}
\newcommand{\cX}{\mathcal{X}}
\newcommand{\dx}{\mathrm{d}x}
\newcommand{\dt}{\mathrm{d}t}
\newcommand{\ds}{\mathrm{d}s}

\numberwithin{equation}{section}

\title[Some one-dimensional elliptic problems with constraints]{Some one-dimensional elliptic problems with constraints}

\author[Jacopo Schino]{Jacopo Schino}
\address[Jacopo Schino]{\newline \indent Faculty of Mathematics, Informatics and Mechanics, University of Warsaw, ul. Banacha 2, 02-097 Warsaw, Poland
\newline \indent and
\newline \indent Neuromed Mediterranean University of Sciences and Technologies, Via dell'Elettronica 1, 86077 Pozzilli (IS), Italy}
\email{\href{mailto:jacopo.schino@unineuromed.it}{jacopo.schino@unineuromed.it}}

\author[Panayotis Smyrnelis]{Panayotis Smyrnelis}
\address[Panayotis Smyrnelis]{\newline \indent Department of Mathematics, University of Athens, 11584 Athens, Greece}
\email{\href{mailto:smpanos@math.uoa.gr}{smpanos@math.uoa.gr}}

\subjclass[2020]{34L40, 35B32, 35J15, 35J35, 35J91}

\keywords{Poly-harmonic Schr\"odinger equations, homoclinic solutions, normalised solutions, least-energy solutions, bifurcation theory, variational methods}

\begin{document}
	\maketitle
	
\begin{abstract}
Given $m \in \mathbb{N} \setminus \{0\}$ and $\rho > 0$, we find solutions $(\lambda,u)$ to the problem
\begin{equation*}
\begin{cases}
\bigl(-\frac{\mathrm{d}^2}{\mathrm{d} x^2}\bigr)^m u + \lambda G'(u) = F'(u)\\
\int_{\mathbb{R}} K(u) \, \mathrm{d}x = \rho
\end{cases}
\end{equation*}
in the following cases: $m=1$ or $2G(s) = K(s) = s^2$. In the former, we follow a bifurcation argument; in the latter, we use variational methods.
\end{abstract}

\section{Introduction and statement of the results}

We study the problem
\begin{equation}\label{eq:main}
\begin{cases}
\displaystyle \Bigl(-\frac{\mathrm{d}^2}{\dx^2}\Bigr)^m u + \lambda G'(u) = F'(u)\\
\displaystyle \int_{\R} K(u) \, \dx = \rho,
\end{cases}
\end{equation}
where $1 \le m \in \mathbb{N}$, $\rho \in (0,+\infty)$ is a prescribed quantity, $\lambda \in \R$ is part of the unknown, and $F$, $G$, and $K$ are suitable functions. 

To explain our motivations, let us start with taking $2G(s) = K(s) = s^2$. In this case, \eqref{eq:main} takes the form
\begin{equation}\label{eq:L2}
\begin{cases}
\displaystyle \Bigl(-\frac{\mathrm{d}^2}{\dx^2}\Bigr)^m u + \lambda u = F'(u)\\
\displaystyle \int_{\R} u^2 \, \dx = \rho.
\end{cases}
\end{equation}
If $F(u) = F(|u|)$, then \eqref{eq:L2} arise when seeking \textit{standing-wave} solutions to the Schr\"odinger-type evolution equation
\begin{equation*}
\mathrm{i} \frac{\partial}{
\partial t} \Psi = \Bigl(-\frac{\partial^2}{\partial x^2}\Bigr)^m \Psi - F'(u),
\end{equation*}
i.e., solutions of the form $\Psi(t,x) = e^{\mathrm{i} \lambda t} u(x)$ with $u\colon \R \to \R$. Then, the $L^2$ constraint is justified because $|\Psi(t,\cdot)| = |u|$ for all $t \in \R$. Solutions to \eqref{eq:L2} are often referred to as \textit{normalised solutions}.

A classical approach to solve \eqref{eq:L2} consists of finding critical points of the functional $J \colon H^m(\R) \to \R$,
\begin{equation}\label{eq:J}
J(u) = \int_{\R} \frac12 |u^{(m)}|^2 - F(u) \, \dx,
\end{equation}
restricted to the set
\begin{equation}\label{eq:S}
\cS := \Set{v \in H^m(\R) | \int_{\R} v^2 \, \dx = \rho}
\end{equation}
under suitable conditions on $F$ that include $F'(s) = \cO(|s|)$ as $s \to 0$.
With this approach, $-\lambda$ is nothing but the Lagrange multiplier arising from the constraint $\cS$. Since minimisers are among the simplest examples of critical points, it makes sense to wonder whether $J|_\cS$ is bounded below: this is determined by the behaviour at infinity of $F(s)$ with respect to $|s|^{2+4m}$ and, sometimes, $\rho$. In particular,
\begin{equation*}
\lim_{|s| \to +\infty} F(s) s^{-(2+4m)} \begin{cases}
\le 0\\
\in (0,+\infty)\\
= +\infty
\end{cases}
\Longrightarrow \quad \inf_\cS J > -\infty \begin{cases}
\text{for all values of } \rho\\
\text{for small values of } \rho\\
\text{for no values of } \rho
\end{cases}
\end{equation*}
(for the sake of the explanation, we assume that $\lim_{|s| \to +\infty} F(s) s^{-(2+4m)}$ exists). These three regimes are known in the literature as mass- (or $L^2$-) subcritical, critical, and supercritical.

It is evident, then, that the number
\begin{equation}\label{eq:ce}
2 + 4m
\end{equation}
($2 + 4m/N$ in dimension $N\ge1$) plays an important role in the geometry of $J|_\cS$, which is why it is called the mass- (or $L^2$-) critical exponent.

When $m=1$, the mass-subcritical case, as well as the mass-critical one with $\rho$ small, were first studied by C.A. Stuart \cite{Stuart} and P.-L. Lions \cite{Lions2}; more recently, they have been dealt with, e.g., in \cite{JL_min,Schino,Shibata}, see also the references therein. In the mass-supercritical regime, instead, the seminal work was carried out by L. Jeanjean \cite{Jeanjean}; lately, the problem was revisited, e.g., in \cite{JL2020CVPD}, see also the references therein. As for the mixed case, it was considered only a few years ago by N. Soave \cite{SoaveJDE} and L. Jeanjean \& S.-S. Lu \cite{Jeanjean_Lu}.

When $m \ge 2$, instead, there is very little work: to our best knowledge, \cite{LuZhang,Phan,ZhangLiuGuan} are the only papers debating the one-dimensional case, and uniquely for $m=2$.

Going back to the article \cite{Lions2}, it is interesting that, for $m=1$, problems more general than \eqref{eq:L2} are considered therein; for example, the author considers \eqref{eq:main} with $pG(s) = K(s) = |s|^p$, $p>1$. At the same time, the recent article \cite{Jean2} introduced a new, non-variational method to look for solutions to \eqref{eq:L2} with $m=1$, which allows the mass-subcritical, -critical, and -supercritical regimes to be dealt with in the same way and where the starting point is the existence of a positive solution to
\begin{equation*}
-u'' + \lambda u = F'(u) \quad \text{in } \R
\end{equation*}
with $\lambda> 0$ fixed and suitable assumptions about $F$. These considerations motivates us to exploit the techniques of \cite{Jean2} and find solutions to \eqref{eq:main} with $G$ and $K$ even more general than in \cite{Lions2}, at least when $m=1$. When $m\ge2$, instead, this new approach does not seem to work because of the lack of a theory about solutions to the differential equation in \eqref{eq:main} with $\lambda$ fixed, and we have to rely on variational methods.

When using an approach inspired from \cite{Jean2}, we consider the following assumptions.

\begin{enumerate}[label=(A\arabic{*}),ref=A\arabic{*}]\setcounter{enumi}{-1}
\item \label{A0} $F,G \in \cC^1([0,+\infty))$, $F(0) = G(0) = F'(0) = G'(0) = 0$, $G'(s) > 0$ for all $s>0$, $\lim_{s\to 0^+}\frac{F(s)}{G(s)} = 0$, and $\lim_{s\to +\infty}\frac{F(s)}{G(s)} = +\infty$. 
\item \label{A1} For all $s > 0$ such that $F(s)>0$ there holds $Z(s) := \big(\frac{F}{G}\big)'(s) > 0$.
\item \label{A2} $K\in \cC([0,\infty))$, $K(0)=0$, $K(s)> 0$ for all $s > 0$, $s \mapsto \frac{K(s)}{\sqrt{G(s)}}$ is integrable in a right-hand neighbourhood of $0$, and $\int_0^{+\infty} \frac{K(s)}{\sqrt{G(s)}} \, \ds = +\infty$.
\end{enumerate}

Moreover, we define
\begin{align*}
\Phi(t) := \left(\int_{0}^{t} \frac{K(s)}{\sqrt{G(s)}} \, \ds\right)^2 \ \text{for } t>0, \quad m_0 := \max\Set{t\geq 0 | F \leq 0 \text{ on } [0,t]},\\
I_F:=\sqrt{2}\int_0^{m_0} \frac{K(s)}{\sqrt{|F(s)|}} \, \ds \in (0,+\infty] \text{ (if $m_0 > 0$),} \qquad \qquad
\end{align*}
and, recalling the definition of $Z$ from (\ref{A1}),
\begin{align*}
& L_0 := \limsup_{s\to 0^+} \frac{Z(s)}{\Phi'(s)} \in [0,+\infty], && \ell_0 := \liminf_{s\to 0^+} \frac{Z(s)}{\Phi'(s)} \in [0,+\infty],\\
& L_\infty := \limsup_{s\to \infty} \frac{Z(s)}{\Phi'(s)} \in [0,+\infty], && \ell_\infty := \liminf_{s\to \infty} \frac{Z(s)}{\Phi'(s)} \in [0,+\infty].
\end{align*}

In this context, our main result reads as follows.

\begin{Th}\label{th:m1}
Let $m=1$ and assume that (\ref{A0})--(\ref{A2}) hold. Then, a solution $(\lambda,u) \in (0,+\infty) \times \cC^2(\R)$ to \eqref{eq:main} exists in each of the following cases:
\begin{itemize}
	\item $m_0=0$, $L_0 < \ell_\infty$, and $\rho \in \left(\frac{\pi}{\sqrt{2  \ell_\infty}}, \frac{\pi}{\sqrt{2 L_0}}\right)$;
	\item $m_0=0$, $L_\infty < \ell_0$, and $\rho \in \left(\frac{\pi}{\sqrt{2  \ell_0}}, \frac{\pi}{\sqrt{2 L_\infty}}\right)$;
	\item $m_0>0$, $I_F > \frac{\pi}{\sqrt{2 \ell_\infty}}$, and $\rho \in \left(\frac{\pi}{\sqrt{2 \ell_\infty}}, I_F\right)$;
	\item $m_0>0$, $F'(m_0) \ne 0$, $I_F < \frac{\pi}{\sqrt{2 L_\infty}}$, and $\rho \in \left(I_F, \frac{\pi}{\sqrt{2 L_\infty}}\right)$.
\end{itemize}
Furthermore, $u' \in L^2(\R)$, $u$ is non-negative and even, and $\lim_{|x| \to +\infty} u(x) = 0$.
\end{Th}

As examples for such $G$ and $K$ we propose
\[
G(s) = \frac1p s^p \quad \text{and} \quad K(s) = s^q, 
\]
with $p>1$ and $q > \max\{p/2-1,0\}$ (observe that this includes the case $p=q=2$). Then, from Theorem \ref{th:m1} we obtain immediately the following outcome.

\begin{Cor}\label{cor:m1}
Let $m=1$, $p>1$, $q > \max\{p/2-1,0\}$, and assume that
\begin{enumerate}[label=(a\arabic{*}),ref=a\arabic{*}]\setcounter{enumi}{-1}
	\item \label{a0} $F \in \cC^1([0,+\infty))$, $\lim_{s\to 0^+} \frac{F(s)}{s^p} = 0$, and $\lim_{s\to +\infty} \frac{F(s)}{s^p} = +\infty$.
	\item \label{a1} For all $s > 0$ such that $F(s)>0$ there holds $F'(s)s - pF(s)>0$.
\end{enumerate}
For $s>0$, define the quantities
\begin{align*}
& K_0 := \limsup_{s\to 0^+} \frac{F'(s)s - pF(s)}{s^{2q+2}}, && k_0 := \liminf_{s\to 0^+} \frac{F'(s)s - pF(s)}{s^{2q+2}},\\
& K_\infty := \limsup_{s\to \infty} \frac{F'(s)s - pF(s)}{s^{2q+2}}, && k_\infty := \liminf_{s\to \infty} \frac{F'(s)s - pF(s)}{s^{2q+2}}.
\end{align*}
Then, a solution $(\lambda,u) \in (0,+\infty) \times \cC^2(\R)$ to \eqref{eq:main} exists in every of the following cases:
\begin{itemize}
	\item $m_0=0$, $K_0 < k_\infty$, and $\rho \in \left(\frac{\pi}{\sqrt{(q-p/2+1)  k_\infty}}, \frac{\pi}{\sqrt{(q-p/2+1) K_0}}\right)$;
	\item $m_0=0$, $K_\infty < k_0$, and $\rho \in \left(\frac{\pi}{\sqrt{(q-p/2+1)  k_0}}, \frac{\pi}{\sqrt{(q-p/2+1) K_\infty}}\right)$;
	\item $m_0>0$, $I_F > \frac{\pi}{\sqrt{(q-p/2+1) k_\infty}}$, and $\rho \in \left(\frac{\pi}{\sqrt{(q-p/2+1) k_\infty}}, I_F\right)$;
	\item $m_0>0$, $F'(m_0) \ne 0$, $I_F < \frac{\pi}{\sqrt{(q-p/2+1) K_\infty}}$, and $\rho \in \left(I_F, \frac{\pi}{\sqrt{(q-p/2+1) K_\infty}}\right)$.
\end{itemize}
Furthermore, $u' \in L^2(\R)$, $u$ is non-negative and even, and $\lim_{|x| \to +\infty} u(x) = 0$.
\end{Cor}

\begin{Rem}\label{rem:JZZ}
Here are some observations about the case $p=q=2$ in Corollary \ref{cor:m1}.
\begin{itemize}
	\item [(i)] The number $2q+2=6$, which appears in the definition of $K_0$, $k_0$, $K_\infty$, and $k_\infty$, is exactly the exponent \eqref{eq:ce} with $m=1$.
	\item [(ii)] At first glance, when $F > 0$ on $(0,+\infty)$, Corollary \ref{cor:m1} is weaker than \cite[Theorem 1.1]{Jean2} in dimension $1$ because of (\ref{a1}); however, such an assumption is needed (cf. Remark \ref{rem:gapJZZ} below), hence \cite[Theorem 6.1]{Jean2}, which \cite[Theorem 1.1]{Jean2} is based on, contains a (small) gap.
	\item [(iii)] Since we can admit sign-changing terms $F$, Corollary \ref{cor:m1} improves the one-dimensional case of \cite[Theorem 1.1]{Jean2}. Additionally, it extends the one-dimensional existence results in \cite{SoaveJDE} to the case of non-linearities more general than the sum of two powers.
\end{itemize}
\end{Rem}

Now, we turn to the case where $m$ is any positive integer. We begin by considering assumptions that describe the mass-subcritical and -critical cases.

\begin{enumerate}[label=(f\arabic{*}),ref=f\arabic{*}]
\setcounter{enumi}{-1}
\item \label{f0} $F \in \cC^1(\R)$ and $F'(s) = \cO(|s|)$ as $s \to 0$.
\item \label{f1} $\lim_{s\to 0} F(s) s^{-2} = 0$.
\item \label{f2} $\sigma := \limsup_{|s| \to +\infty} F(s)s^{-(2+4m)} < +\infty$.
\item \label{f3} $\lim_{s\to 0} F(s)s^{-(2+4m)} = +\infty$.
\end{enumerate}

We recall from \eqref{eq:J} and \eqref{eq:S} the definitions of $J$ and $\cS$ and introduce the set
\begin{equation*}
\cD := \Set{v \in H^m(\R) | \int_{\R} v^2 \, \dx \le \rho},
\end{equation*}
which was first used in \cite{BiegMed} and then exploited, e.g., in \cite{BdS,BMS,BMiS,CLY,Liu_Zhao_CVPD,MedSc,MSlog,Schino,LRZ}. Moreover, let us recall the Gagliardo--Niremberg inequality \cite{Gagliardo1,Gagliardo2,Nirenberg}, here expressed in the one-dimensional case: for every $p > 2$ there exists $C_p > 0$ such that for all $v \in H^m(\R)$ there holds
\begin{equation}\label{eq:GN}
|v|_p \le C_p |v^{(m)}|_2^{\delta_p} |v|_2^{1-\delta_p},
\end{equation}
and $C_p$ is sharp, where $\delta_p = (1/2 - 1/p) / m$ and $|\cdot|_q$ denotes the norm in $L^q(\R)$, $q \in [1,+\infty]$.

Our existence result in this regime is the following.

\begin{Th}\label{th:mainSUB}
If (\ref{f0})--(\ref{f3}) are satisfied and $2\sigma C_{2+4m}^{2+4m} \rho^{2m} < 1$, then there exist $u \in \cS$ and $\lambda > 0$ such that $J(u) = \min_{\cD} J < 0$ and $(\lambda,u)$ is a solution to \eqref{eq:L2}.
\end{Th}

The proof of Theorem \ref{th:mainSUB} follows verbatim that of \cite[Theorem 1.1]{Schino} once proved that every solution $u \in H^m(\R)$ to the differential equation in \eqref{eq:L2}
with $\lambda \in \R$ fixed satisfies the Poho\v{z}aev identity if (\ref{f0}) holds, which is done in Proposition \ref{pr:Po} below. For this reason, we omit it.

Next, we move to the mass-supercritical case. Let us define $H(s) := F'(s)s - 2F(s)$ for $s \in \R$. We assume the following conditions. 

\begin{enumerate}[label=(F\arabic{*}),ref=F\arabic{*}]
\setcounter{enumi}{-1}
\item \label{F0} $F,H \in \cC^1(\R)$ and $|F'(s)| + |H'(s)| = \cO(|s|)$ as $s \to 0$.
\item \label{F1} $\eta := \limsup_{s \to 0} H(s)s^{-(2+4m)} < +\infty$. 
\item \label{F2} $\lim_{|s| \to +\infty} F(s) s^{-(2+4m)} = +\infty$
\item \label{F3} $(2+4m) H(s) \le H'(s) s$ for all $s \in \R$.
\item \label{F4} $0\le 4m F(s) \le H(s)$ for all $s \in \R$.
\end{enumerate}

Since $J|_\cS$ is unbounded below under (\ref{F0})--(\ref{F4}), cf. the proof of Lemma \ref{l:fiber} below, we follow the approach in \cite{BiegMed}, where -- moreover -- examples of such functions $F$ can be found.

Let us recall the Nehari and Poho\v{z}aev identities associated with \eqref{eq:L2}, i.e.,
\begin{align*}
\int_{\R} |u^{(m)}|^2 + \lambda u^2 \, \dx & = \int_{\R} F'(u) u \, \dx\\
\int_{\R} (1-2m) |u^{(m)}|^2 + \lambda u^2 \, \dx & = 2 \int_{\R} F(u) \, \dx
\end{align*}
respectively. Then, every $u \in H^m(\R) \setminus \{0\}$ that solves the differential equation in \eqref{eq:L2} for some $\lambda \in \R$ belongs to the set
\begin{equation*}
\cM := \Set{v \in H^m(\R) \setminus \{0\} | \int_{\R} |v^{(m)}|^2 \, \dx = \frac{1}{2m} \int_{\R} H(v) \, \dx}.
\end{equation*}

If $H \in \cC^1(\R)$, $H'(s) = \cO(|s|)$ as $s \to 0$, and $H(\xi_0) > 0$ for some $\xi_0 \ne 0$, one easily proves that $\cM$ is a manifold of class $\cC^1$ and co-dimension $1$; see, e.g., \cite[Lemma 4.1]{BMS}.

We consider the following condition, which will be paired with (\ref{F3}):
\begin{equation}\label{eq:strict}
\int_\R H'(u) u - (2+4m) H(u) \, \dx>0 \quad \forall u \in H^m(\R) \setminus  \{0\}.
\end{equation}
Note that \eqref{eq:strict} is satisfied if and only if  $H(s)s - (2+4m) H(s)\ge 0$ holds for every $s \in \R$ and the strict inequality holds along two sequences $s'_n \to 0^+$ and $s''_n \to 0^-$ (cf. \cite[Lemma 2.1]{BiegMed}).

Finally, we introduce the condition
\begin{equation}\label{eq:rho}
\eta C_{2+4m}^{2+4m} \rho^{2m} < 2m,
\end{equation}
which allows us to deal with a non-linearity with mass-critical growth at the origin.

Our existence result in this regime is the following.

\begin{Th}\label{th:mainSUP}
If (\ref{F0})--(\ref{F4}) and \eqref{eq:rho} hold, then there exists $u \in \cD \cap \cM$ such that $J(u) = \min_{\cD \cap \cM} J > 0$. If, moreover, \eqref{eq:strict} holds, then there exists $\lambda > 0$ such that $(\lambda,u)$ is solution to \eqref{eq:main} -- in particular, $u \in \cS$.
\end{Th}

The proof of Theorem \ref{th:mainSUP} is modelled on that of \cite[Theorem 3.3]{BMS}; nonetheless, since the one-dimensional setting requires some modifications, we provide it in Section \ref{s:vm} for the reader's convenience. Theorem \ref{th:m1}, instead, is proved in Section \ref{s:m1}.

\section{The global-branch approach}\label{s:m1}

Throughout this section, $m=1$, and we assume (\ref{A0})--(\ref{A2}).

In view of (\ref{A0}), for every $\lambda>0$ there exists $m_\lambda \in (0,+\infty)$ such that $W_\lambda:=\lambda G-F$ is positive on $(0,m_\lambda)$ and $W_\lambda(0) = W'_\lambda(0) = W_\lambda(m_\lambda) = 0$. On the other hand, (\ref{A1}) implies that $W'_\lambda(m_\lambda)=-Z(m_\lambda)G(m_\lambda)<0.$ Consequently, for every $\lambda>0$, there exists a homoclinic orbit $u_\lambda$ corresponding to $W_\lambda$ (cf. for instance \cite[Theorem 5]{BerLionsI} or \cite[Theorem 5.4]{Smy}\footnote{Both these theorems require additional regularity for $F$ and $G$ (at least $\cC_\textup{loc}^{1,1}$); however, similar arguments as in \cite{Smy} can be repeated with minor modifications when $F$ and $G$ are merely $\cC^1$.}) satisfying the following properties:
\begin{itemize}
	\item  $u_\lambda \in\cC^2(\R)$ is a non-negative and even solution of $-u'' + \lambda G'(u) = F'(u)$.
	\item $\lim_{|x|\to\infty} u_\lambda(x) = 0$ and $m_\lambda = u_\lambda(0) = \max_\R u_\lambda$.
	\item $|u'_\lambda(x)|^2 = 2W_\lambda(u_\lambda(x))$ for all $x\in\R$ (equipartition relation).
	\item Setting $T_\lambda:=\int_0^{m_\lambda}\frac{\mathrm{d} u}{\sqrt{2W_\lambda(u)}}\in (0,+\infty]$\footnote{Observe that $T_\lambda = +\infty$ if $W_\lambda \in \cC^{1,1}([0,\varepsilon])$ for some $\varepsilon > 0$.}, $u_\lambda$ is increasing on $(-T_\lambda,0)$ and decreasing on $(0,T_\lambda)$, while $u_\lambda(x)=0$, if $|x|\geq T_\lambda$.
	\item $u'_\lambda\in L^2(\R)$.

\end{itemize}

\begin{Rem}\label{rem:gapJZZ}
We point out that the assumption $F'(t)>0$ for all $t>0$ is not sufficient to ensure, for every $\lambda>0$, the existence of the homoclinic orbit $u_\lambda$, since  we may have $W'_\lambda(m_\lambda)= 0$ for some  $\lambda>0$. For example, when $G(s) = s^2/2$, taking
\begin{equation*}
F(s) =
\begin{cases}
\frac12 s^2 + \cos(s) - 1 & \text{if } s \in [0,2\pi]\\
\frac12 s^2 + (s-2\pi)^p & \text{if } s > 2\pi
\end{cases}
\end{equation*}
with $p>2$ we see that $m_1 = 2\pi$ and $W_1'(m_1) = 0$, thus no non-trivial non-negative solutions to $-u'' + u = F'(u)$ that vanish at infinity exist in view of \cite[Theorem 5]{BerLionsI}.
\end{Rem}

Recalling that $m_0=\max\{t\geq 0: F\leq 0 \text{ on } [0,t]\}$, we have $F(m_0)=0$, and in view of (\ref{A1}), $F$ and $F'$ are positive on $(m_0,+\infty)$. In addition, the following properties hold.

\begin{Lem}\label{lem:ml}
The function $(0,+\infty) \ni \lambda \mapsto m_\lambda \in (0,+\infty)$ is increasing, continuous, $\lim_{\lambda\to+\infty} m_\lambda = +\infty$, and $\lim_{\lambda\to 0^+} m_\lambda = m_0$.
\end{Lem}
\begin{proof}
Observe that (\ref{A0}) implies that $G(s) > 0$ if $s > 0$. Let $0 < \lambda < \Lambda < +\infty$ and $s \in (0,m_\lambda]$; then,
\[
W_\Lambda(s) = \Lambda G(s) - F(s) > \lambda G(s) - F(s) = W_\lambda(s) \ge 0,
\]
which shows that $m_\Lambda > m_\lambda$.\\
The continuity follows from the implicit function theorem applied to the function $(\lambda,s) \mapsto W_\lambda(s)$ because, as remarked above, $W'_\lambda(m_\lambda) < 0$.\\
Next, observe that the two limits in the statement exist in virtue of the monotonicity of $\lambda \mapsto m_\lambda$. Assume by contradiction that $\lim_{\lambda\to+\infty} m_\lambda =: M < +\infty$. Then
\[
0 = \lim_{\lambda\to+\infty} W_\lambda(m_\lambda) = \lim_{\lambda\to+\infty} \lambda G(m_\lambda) - F(m_\lambda) = +\infty,
\]
which is impossible.\\
Finally, observe that the definition of $m_0$ and the fact that $F(m_\lambda) = m_\lambda G(m_\lambda) > 0$ imply that $m_0 < m_\lambda$. Denoting $M := \lim_{\lambda\to 0^+} m_\lambda \ge m_0$, we have
\[
0 = \lim_{\lambda\to 0^+} W_\lambda(m_\lambda) = \lim_{\lambda\to 0^+} \lambda G(m_\lambda) - F(m_\lambda) = - F(M),
\]
and the statement follows from the fact that $F$ is positive on $(m_0,+\infty)$.
\end{proof}

\begin{Lem}\label{lem:rl}
For every $\lambda > 0$ there holds
\[
\rho_\lambda := \int_\R K( u_\lambda(x)) \, \dx = \sqrt{2} \int_0^{m_\lambda} \frac{K(u)}{\sqrt{G(u)}} \left(\frac{F(m_\lambda)}{G(m_\lambda)} - \frac{F(u)}{G(u)}\right)^{-1/2} \, \mathrm{d}u < +\infty.
\]
Furthermore, $\lambda \mapsto \rho_\lambda$ is continuous.
\end{Lem}
\begin{proof}
Observe that $u_\lambda|_{(-T_\lambda,0)} \colon (-T_\lambda,0) \to (0,m_\lambda)$ is a diffeomorphism and, from the equipartition relation, $u_\lambda'(x) = \sqrt{2 W_\lambda(u_\lambda(x))}$. Consequently,
\begin{equation}\label{eq:rl1}
\rho_\lambda = \int_\R K(u_\lambda(x)) \, \dx = 2\int_{-T_\lambda}^0 K(u_\lambda(x)) \, \dx = \sqrt{2} \int_0^{m_\lambda} \frac{K(u)}{\sqrt{W_\lambda(u)}} \, \mathrm{d}u.
\end{equation}
Since, from (\ref{A0}) and (\ref{A2}), $F(s) / G(s) \to 0$ as $s \to 0^+$ and $K / \sqrt{G}$ is integrable in a right-hand neighbourhood of $0$, we have that $K / \sqrt{W_\lambda}$ is integrable in a right-hand neighbourhood of $0$ as well. This, together with the property that $W_\lambda'(m_\lambda) \ne 0$, yields that $\rho_\lambda < +\infty$. Moreover,  the continuity of $\lambda \mapsto \rho_\lambda$ follows from the one of $\lambda \mapsto m_\lambda$ (Lemma \ref{lem:ml}) and the facts above. Finally,
\begin{equation}\label{eq:rl2}
\begin{split}
\int_0^{m_\lambda} \frac{K(u)}{\sqrt{W_\lambda(u)}} \, \mathrm{d}u & = \int_0^{m_\lambda} \frac{K(u)}{\sqrt{G(u)}} \left(\lambda - \frac{F(u)}{G(u)}\right)^{-1/2} \, \mathrm{d}u\\
& = \int_0^{m_\lambda} \frac{K(u)}{\sqrt{G(u)}} \left(\frac{F(m_\lambda)}{G(m_\lambda)} - \frac{F(u)}{G(u)}\right)^{-1/2} \, \mathrm{d}u,
\end{split}
\end{equation}
and the conclusion follows from \eqref{eq:rl1}--\eqref{eq:rl2}.
\end{proof}

Next,  we determine the asymptotic behaviour of $\rho_\lambda$ as $\lambda \to 0^+$ and $\lambda \to +\infty$.

\begin{Lem}\label{asympb}
Assume (\ref{A0})--(\ref{A2}). With the convention that $\frac{1}{0} = +\infty$ and $\frac{1}{+\infty} = 0$, the following holds.
\begin{enumerate}[label=(\arabic{*}),ref=\arabic{*}]
	\item \label{1} If $m_0=0$, then $\liminf_{\lambda\to 0^+} \rho_\lambda \geq \frac{\pi}{\sqrt{2L_0}}$ and $\limsup_{\lambda\to 0^+} \rho_\lambda \leq \frac{\pi}{\sqrt{2\ell_0}}$.
	\item \label{2} If $m_0>0$, then $\liminf_{\lambda\to 0^+} \rho_\lambda \geq  I_F$. If, in addition, $F'(m_0) \ne 0$,
	then $\lim_{\lambda\to 0^+} \rho_\lambda = I_F$.
	\item \label{3} $\liminf_{\lambda\to+\infty} \rho_\lambda \geq \frac{\pi}{\sqrt{2L_\infty}}$ and $\limsup_{\lambda\to+\infty} \rho_\lambda \leq  \frac{\pi}{\sqrt{2\ell_\infty}}$.
\end{enumerate}	
\end{Lem}
\begin{proof}
(\ref{1}) Given $L \in (L_0,+\infty)$, there exists $\delta_L>0$ such that $Z(s) \leq L \Phi'(s)$ for all $s \in (0,\delta_L)$. In addition, from Lemma \ref{lem:ml}, $m_\lambda<\delta_L$ holds for $\lambda\in (0,\Lambda_L)$, with $\Lambda_L>0$ small enough. Therefore, we have for $\lambda\in (0,\Lambda_L)$
$$
\frac{F(m_\lambda)}{G(m_\lambda)}-\frac{F(u)}{G(u)}\leq L \bigl(\Phi(m_\lambda)-\Phi(u)\bigr) \quad \forall u\in[0,m_\lambda]
$$
and, using Lemma \ref{lem:rl} and the change of variable $t = \sqrt{\frac{\Phi(u)}{\Phi(m_\lambda)}}$,
\begin{equation*}
\rho_\lambda \geq \sqrt{\frac{2}{L}} \int_0^{m_\lambda} \frac{K(u)(G(u))^{-1/2}}{(\Phi(m_\lambda)-\Phi(u))^{1/2}} \, \mathrm{d}u = \sqrt{\frac{2}{L}} \int_0^{1} \frac{1}{(1-t^2)^{1/2}} \, \dt = \frac{\pi}{\sqrt{2L}},
\end{equation*}
which proves the first part of (\ref{1}). The second one is proved similarly.
	
(\ref{2}) From \eqref{eq:rl1} and the monotone convergence theorem,
$$
\liminf_{\lambda\to 0^+} \rho_\lambda \geq \lim_{\lambda\to 0^+} \sqrt{2} \int_0^{m_0} \frac{K(u)}{\sqrt{W_\lambda(u)}} \, \mathrm{d}u = \sqrt{2} \int_0^{m_0} \frac{K(u)} {\sqrt{|F(u)|}} \, \mathrm{d}u = I_F.
$$
Next, if $F'(m_0) \ne 0$, which implies $Z(m_0) \ne 0$, then -- cf. \eqref{eq:rl2}
\begin{align*}
\int_{m_0}^{m_\lambda} \frac{K(u)}{\sqrt{W_\lambda(u)}} \, \mathrm{d}u & = \int_{m_0}^{m_\lambda} \frac{K(u)(G(u))^{-1/2}}{\sqrt{\frac{F(m_\lambda)}{G(m_\lambda)} - \frac{F(u)}{G(u)}}} \, \mathrm{d}u = \int_{m_0}^{m_\lambda} \frac{K(u)(G(u))^{-1/2}}{\sqrt{\int_u^{m_\lambda} Z(s) \, \ds}} \, \mathrm{d}u\\
& = \cO(\sqrt{m_\lambda-m_0})\to 0 \quad \text{as } \lambda \to 0^+,
\end{align*}
therefore $\lim_{\lambda\to 0^+}\rho_\lambda = I_F$.

(\ref{3}) Let us begin with the first part. Given $L \in (L_\infty,+\infty)$, there exists $S_L > m_0$ such that $Z(s)\leq L \Phi'(s)$ for all $s\geq S_L$. In addition, from Lemma \ref{lem:ml}, $m_\lambda>S_L$ holds for $\lambda\in (\Lambda_L,\infty)$, with $\Lambda_L>0$ large enough. Therefore, we have for $\lambda\in (\Lambda_L,+\infty)$
\begin{equation*}
\frac{F(m_\lambda)}{G(m_\lambda)}-\frac{F(u)}{G(u)} \leq L \bigl(\Phi(m_\lambda)-\Phi(u)\bigr) \quad \text{for all } u \in [S_L,m_\lambda].
\end{equation*}
This, together with (\ref{A2}) and Lemmas \ref{lem:ml} and \ref{lem:rl}, implies
\begin{align*}
\rho_\lambda & \ge \sqrt{2} \int_{S_L}^{m_\lambda} \frac{K(u)}{\sqrt{G(u)}} \left(\frac{F(m_\lambda)}{G(m_\lambda)} - \frac{F(u)}{G(u)}\right)^{-1/2} \, \mathrm{d}u\\
& \ge \sqrt{\frac2L} \int_{S_L}^{m_\lambda} \frac{K(u)(G(u))^{-1/2}}{(\Phi(m_\lambda) - \Phi(u))^{1/2}} \, \mathrm{d}u = \sqrt{\frac2L} \int_{\sqrt{\Phi(S_L) / \Phi(m_\lambda)}}^1 \frac{1}{(1-t^2)^{1/2}} \, \dt\\
& \to \sqrt{\frac2L} \int_0^1 \frac{1}{(1-t^2)^{1/2}} \, \dt = \frac{\pi}{\sqrt{2L}} \quad \text{as } \lambda \to +\infty.
\end{align*}

Now we move to the second part. Given $L \in (0,\ell_\infty)$, there exists $S_L > m_0$ such that $Z(s)\geq L \Phi'(s)$ for all $s\geq S_L$. Additionally, since $F/G$ is non-positive on $[0,m_0]$ and, from (\ref{A1}), positive and increasing on $(m_0,+\infty)$, there holds
$$
\frac{F(S_L)}{G(S_L)} \ge \frac{F(u)}{G(u)} \quad \text{for all } u\in(0,S_L].
$$
In addition, from Lemma \ref{lem:ml}, $m_\lambda>S_L$ holds for $\lambda\in (\Lambda_L,\infty)$, with $\Lambda_L>0$ large enough. Therefore, we have for $\lambda\in (\Lambda_L,+\infty)$
\begin{align*}
\frac{F(m_\lambda)}{G(m_\lambda)}-\frac{F(u)}{G(u)} & \geq L \bigl(\Phi(m_\lambda)-\Phi(u)\bigr) \quad \text{for all } u \in [S_L,m_\lambda],\\
\frac{F(m_\lambda)}{G(m_\lambda)}-\frac{F(u)}{G(u)} & \geq L \bigl(\Phi(m_\lambda)-\Phi(S_L)\bigr) \quad \text{for all } u\in (0,S_L],
\end{align*}
This, together with Lemma \ref{lem:rl}, implies
$$
\rho_\lambda \leq \sqrt{\frac{2}{L}} \bigl(I_1(\lambda)+I_2(\lambda)\bigl) \quad \text{for all } \lambda \in (\Lambda_L,+\infty),
$$
where, thanks to (\ref{A2}),
\begin{align*}
I_1(\lambda) & = \int_0^{S_L}   \frac{K(u)(G(u))^{-1/2}}{(\Phi(m_\lambda) - \Phi(S_L))^{1/2}} \, \mathrm{d}u \to 0 \quad \text{as } \lambda \to +\infty,\\
I_2(\lambda) & = \int_{S_L}^{m_\lambda} \frac{K(u)(G(u))^{-1/2}}{(\Phi(m_\lambda) - \Phi(u))^{1/2}} \, \mathrm{d}u \leq \int_{0}^{m_\lambda} \frac{K(u)(G(u))^{-1/2}}{(\Phi(m_\lambda) - \Phi(u))^{1/2}} \, \mathrm{d}u = \frac{\pi}{2},
\end{align*}
which proves the second part of (\ref{3}).
\end{proof}

\begin{proof}[Proof of Theorem \ref{th:m1}]
It follows from Lemmas \ref{lem:rl} and \ref{asympb}.
\end{proof}

\begin{Rem}
Similarly, non-existence results for problem \eqref{eq:main} can be established in view of Lemma \ref{asympb} and the continuity of the function $\lambda \mapsto \rho_\lambda$.
\end{Rem}

\section{The poly-harmonic case}\label{s:vm}

We begin by proving that every solution to the differential equation in \eqref{eq:L2} satisfies the Poho\v{z}aev identity. In fact, we can state a more general result.

\begin{Prop}[Poho\v{z}aev identity]\label{pr:Po}
Let $N \in [1,2m-1]$ be an integer, $g \in \cC(\R)$ such that $g(s) = \cO(|u|)$ as $s \to 0$, and define $G(s) := \int_{0}^{s} g(t) \, \dt$. If $u \in H^m(\R)$ is a weak solution to
\begin{equation}\label{eq:eq}
(-\Delta)^m u = g(u) \quad \text{in } \R^N,
\end{equation}
then
\begin{equation*}
(N-2m) \int_{\R^N} |\nabla^m u|^2 \, \dx = 2N \int_{\R^N} G(u) \, \dx.,
\end{equation*}
where
\begin{align*}
\nabla^m u &:=
\begin{cases}
\Delta^{m/2} u & \text{if } m  \text{ is even,}\\
\nabla \Delta^{(m-1)/2} u & \text{if } m \text{ is odd.}
\end{cases}
\end{align*}
\end{Prop}
\begin{proof}
Since, from the Sobolev embedding, $|u|_\infty < +\infty$, we have that $|g(u)| \lesssim |u|$ a.e. in $\R^N$, hence $g \circ u \in L^2(\R^N) \cap L^\infty(\R^N)$. Then, from \cite[Lemma 3.1]{Siemianowski}, $u \in W^{2m,p}_\textup{loc}(\R^N)$ for every $p \in [1,+\infty)$. Now, one can follow \cite[Proof of Proposition 2.5]{BMS}. We sketch the proof for the reader's convenience.\\
For every $n \ge 1$, let $\psi_n \in \cC^1_0(\R^N)$ radially symmetric such that $0 \le \psi_n \le 1$, $\psi_n(x) = 1$ for every $|x| \le n$, $\psi_n(x) = 0$ for every $|x| \ge 2n$, and $|x||\nabla \psi_n(x)| \lesssim 1$ for every $x \in \R^N$.\\
Next, observe that the following identities hold true:
\begin{align*}
g(u) (\nabla u \cdot x) \psi_n  = & \; \nabla \cdot \bigl( \psi_n G(u) x \bigr) - N \psi_n G(u) - G(u) \nabla \psi_n \cdot x,\\
\Delta^{2k+1} u (\nabla u \cdot x) \psi_n = & \; \nabla \cdot \biggl[\biggl(\Delta^k(x \cdot \nabla u) \nabla \Delta^k u - \frac{|\nabla \Delta^k u|^2}{2} x\\
& - \sum_{j=0}^{k-1} \Delta^{2k-j} u \nabla \Delta^j (\nabla u \cdot x) + \sum_{j=0}^{k-1} \Delta^j(\nabla u \cdot x) \nabla \Delta^{2k-j}u\biggr) \psi_n\biggr]\\
& + \frac{N-4k-2}{2} |\nabla \Delta^k u|^2 \psi_n - \biggl(\Delta^k(\nabla u \cdot x) \nabla \Delta^k u - \frac{|\nabla \Delta^k u|^2}{2} x\\
& - \sum_{j=0}^{k-1} \Delta^{2k-j} u \nabla \Delta^j (\nabla u \cdot x) + \sum_{j=0}^{k-1} \Delta^j(\nabla u \cdot x) \nabla \Delta^{2k-j}u\biggr) \cdot \nabla \psi_n,\\
\Delta^{2k}u (\nabla u \cdot x) \psi_n = & \; \nabla \cdot\biggl[\biggl(\frac12 (\Delta^ku)^2 x + (\nabla u \cdot x) \nabla \Delta^{2k-1}u\\
& + \sum_{j=0}^{k-2} \Delta^{j+1}(\nabla u \cdot x) \nabla \Delta^{2k-j-2}u\\
& - \sum_{j=0}^{k-1} \Delta^{2k-j-1}u \nabla \Delta^j(\nabla u \cdot x)\biggr) \psi_n\biggr] + \frac{4k-N}{2} (\Delta^ku)^2 \psi_n\\
& - \biggl(\frac12 (\Delta^ku)^2 x + (\nabla u \cdot x) \nabla \Delta^{2k-1}u\\
& + \sum_{j=0}^{k-2} \Delta^{j+1}(\nabla u \cdot x) \nabla \Delta^{2k-j-2}u\\
& - \sum_{j=0}^{k-1} \Delta^{2k-j-1}u \nabla \Delta^j(\nabla u \cdot x)\biggr) \cdot \nabla \psi_n.
\end{align*}
Multiplying both sides of \eqref{eq:eq} by $\psi_n \nabla u \cdot x$, using the identities above, and integrating over $\R^N$, we obtain
\begin{equation}\label{eq:eqn}
\begin{split}
0 & = \int_{\R^N} \left( -(-\Delta)^m u + g(u) \right) \psi_n \nabla u \cdot x \, \dx\\
& = \int_{\R^N} \frac12 |\nabla^m u|^2 \nabla \psi_n \cdot x + \cX \cdot \nabla \psi_n + \frac{N-2m}{2} \psi_n |\nabla^m u|^2 - N \psi_n G(u)\\
& \qquad - G(u) \nabla \psi_n \cdot x + \nabla \cdot \left[ \psi_n \left( - \cX - \frac12 |\nabla^m u|^2 x + G(u)x \right) \right] \, \dx,
\end{split}
\end{equation}
where
\begin{equation*}
\cX :=
\begin{cases}
\displaystyle - \Delta^k(\nabla u \cdot x) \nabla \Delta^k u + \sum_{j=0}^{k-1} \Delta^{2k-j} u \nabla \Delta^j (\nabla u \cdot x) - \sum_{j=0}^{k-1} \Delta^j(\nabla u \cdot x) \nabla \Delta^{2k-j}u &\\
\displaystyle \nabla u \cdot x \nabla \Delta^{2k-1} u + \sum_{j=0}^{k-2} \Delta^{j+1}(\nabla u \cdot x) \nabla \Delta^{2k-j-2}u - \sum_{j=0}^{k-1} \Delta^{2k-j-1}u \nabla \Delta^j(\nabla u \cdot x)
\end{cases}
\end{equation*}
if $m=2k+1$ or $m=2k$ respectively.\\
Finally, from the properties of $\psi_n$ and the dominated convergence theorem, we conclude the proof letting $n \to +\infty$ in \eqref{eq:eqn}.
\end{proof}

For $u \in H^m(\R) \setminus \{0\}$ and $s>0$, let us define $s \star u := \sqrt{s} \, u(s \cdot)$ and $\varphi_u(s) := J(s \star u)$. Note that $|s \star u|_2 = |u|_2$ and that $s \star u \in \cM$ if and only if $\varphi_u'(s) = 0$.

\begin{Lem}\label{l:fiber}
Assume (\ref{F0})--(\ref{F4}) and \eqref{eq:rho} hold. For every $u \in H^m(\R)$ such that
\begin{equation}\label{eq:fiber}
\eta |u|_{2+4m}^{2+4m} < 2m |u^{(m)}|_2^2
\end{equation}
there exist $0 < a \le b < +\infty$ such that $\varphi_u$ is increasing on $(0,a)$, decreasing on $(b,+\infty)$, and $\varphi_u \equiv \max \varphi_u$ on $[a,b]$. If, moreover, \eqref{eq:strict} is satisfied, then $a=b$.
\end{Lem}
Notice that every $u \in \cD \cap H^m(\R) \setminus \{0\}$ satisfies \eqref{eq:fiber} if \eqref{eq:rho} holds.
\begin{proof}
Let $u \in H^m(\R)$ as in the assumptions. From (\ref{F1}), (\ref{F4}), and the continuity of $F$, there exists $c > 0$ such that
\begin{equation*}
|F(t)| \le c t^{2+4m} \quad \forall t  \in \overline{B(0,|u|_\infty)},
\end{equation*}
whence
\begin{equation*}
\varphi_u(s) = \frac{s^{2m}}{2} |u^{(m)}|_2^2 - \int_{\R} \frac{F(s^{1/2}u)}{s} \, \dx \to 0
\end{equation*}
as $s \to 0^+$. Moreover,
\[
\frac{\varphi_u(s)}{s^{2m}} = \frac{|u^{(m)}|_2^2}{2} - \int_{\R} \frac{F(s^{1/2}u)}{s^{1+2m}} \, \dx ,
\]
with $\lim_{s \to +\infty} \int_{\R} F(s^{1/2}u) / s^{1+2m} \, \dx = +\infty$ from (\ref{F2}) and Fatou's Lemma. This proves that $\lim_{s \to +\infty} \varphi_u(s) = -\infty$.

Now, fix $\varepsilon > 0$ such that $(\eta + 2m\varepsilon) |u|_{2+4m}^{2+4m} < 2m |u^{(m)}|_2^2$. From (\ref{F1}), (\ref{F4}), and the continuity of $F$, there exists $C = C(\varepsilon,|u|_\infty) > 0$ such that
\[
F(t) \le \left(\frac{\eta}{4m} + \varepsilon\right) t^{2+4m} + C t^{4+4m} \quad \forall t \in \overline{B(0,|u|_\infty)},
\]
whence, using also \eqref{eq:fiber},
\begin{align*}
\varphi_u(s) & \ge \frac{s^{2m}}{2} |u^{(m)}|_2^2 - \frac1s \left[\left(\frac{\eta}{4m} + \varepsilon\right) |s^{1/2} u|_{2+4m}^{2+4m} + C |s^{1/2} u|_{4+4m}^{4+4m}\right]\\
& = \frac{s^{2m}}{2} \left[|u^{(m)}|_2^2 - \left(\frac{\eta}{2m} + \varepsilon\right) |u|_{2+4m}^{2+4m}\right] - C |u|_{4+4m}^{4+4m} s^{2m+1},
\end{align*}
which proves that $\varphi_u(s) > 0$ if $s \ll 1$. So far we have proved that $\varphi_u$ attains its positive maximum. To conclude, note that
\[
\varphi_u'(s) = m s^{2m-1} \left(|u^{(m)}|_2^2 - \frac{1}{2m} \int_{\R} \frac{H(s^{1/2} u)}{s^{1+2m}} \, \dx\right),
\]
where the function
\[
s \mapsto \int_{\R} \frac{H(s^{1/2} u)}{s^{1+2m}} \, \dx
\]
is non-decreasing from (\ref{F0}) and (\ref{F3}), and increasing if \eqref{eq:strict} holds. 
\end{proof}

\begin{Rem}\label{r:nonempty}
If  (\ref{F0})--(\ref{F4}) and \eqref{eq:rho} hold, then given $u \in \cS \cap H^m(\R) \setminus \{0\}$, it follows from the proof of Lemma \ref{l:fiber} that  $\varphi_u'(s)=0$ for some $s>0$. Thus, we have $s \star u \in \cS\cap\cM$, and $\cS\cap \cM\neq\emptyset$.
\end{Rem}

\begin{Rem}\label{r:empty}
If (\ref{F0}), (\ref{F3}), and \eqref{eq:strict} hold, then $\Set{u \in \cM | \varphi_u''(1) = 0} = \emptyset$. As a matter of fact, from $\varphi_u'(1) = \varphi_u''(1) = 0$ we obtain
\[
\int_{\R}(2+4m) H(u) - H'(u) u \, \dx = 0,
\]
which contradicts \eqref{eq:strict} because $u \ne 0$.
\end{Rem}

\begin{Lem}\label{l:grad}
If (\ref{F0}), (\ref{F1}), and \eqref{eq:rho} hold, then $\inf_{u \in \cD \cap \cM} |u^{(m)}|_2 > 0$.
\end{Lem}
\begin{proof}
Let $u \in \cD \cap \cM$. If $\|u\|_{H^m} > \sqrt{2\rho}$, then
\[
|u^{(m)}|_2^2 = \|u\|_{H^m}^2 - |u|_2^2 > 2\rho - \rho = \rho,
\]
hence we assume that $\|u\|_{H^m} \le \sqrt{2\rho}$. This implies that $|u|_\infty \le c \sqrt{2\rho}$, where $c>0$ is the best constant in the embedding $H^m(\R) \hookrightarrow L^\infty(\R)$. Let $\varepsilon > 0$:
from (\ref{F0}) and (\ref{F1}) there exists $C = C(\varepsilon,\rho) > 0$ such that $H(t) \le (\eta + \varepsilon)t^{2+4m} + C t^{4+4m}$ for every $t \in \overline{B(0,c \sqrt{2\rho})}$. Consequently, from \eqref{eq:GN} and the fact that $u \in \cM$,
\begin{align*}
2m |u^{(m)}|_2^2 & = \int_{\R} H(u) \, \dx \le (\eta+ \varepsilon) |u|_{2+4m}^{2+4m} + C |u|_{4+4m}^{4+4m}\\
& \le (\eta + \varepsilon) C_{2+4m}^{2+4m} |u^{(m)}|_2^2 |u|_2^{4m} + C C_{4+4m}^{4+4m} |u^{(m)}|_2^{2+1/m} |u|_2^{2+4m-1/m}\\
& \le (\eta+ \varepsilon) C_{2+4m}^{2+4m} \rho^{2m} |u^{(m)}|_2^2 + C C_{4+4m}^{4+4m} \rho^{1+2m-1/(2m)} |u^{(m)}|_2^{2+1/m},
\end{align*}
and we conclude taking $\varepsilon > 0$ sufficiently small in view of \eqref{eq:rho}.
\end{proof}

\begin{Lem}\label{l:positive}
If (\ref{F0})--(\ref{F4}) and \eqref{eq:rho} hold, then $\inf_{\cD \cap \cM} J > 0$.
\end{Lem}
\begin{proof}
We begin by showing the existence of $\delta > 0$ such that
\begin{equation*}
\left(2m - \eta C_{2+4m}^{2+4m} \rho^{2m}\right) |u^{(m)}|_2^2 \le 8m J(u)
\end{equation*}
for all $u \in \cD \cap H^m(\R)$ with $|u^{(m)}|_2 \le \delta$. We can assume $\delta \le 1$, hence there exists $c_\rho > 0$ such that $|u|_\infty \le c_\rho$ for all $u$ as above. From (\ref{F0}), (\ref{F1}), and (\ref{F4}), for every $\varepsilon > 0$ there exists $C = C(\varepsilon,\rho) > 0$ such that
\[
F(t) \le \left(\frac{\eta}{4m} + \varepsilon\right)t^{2+4m} + C t^{4+4m} \quad \forall t \in \overline{B(0,c_\rho)}.
\]
This and \eqref{eq:GN} yield
\begin{align*}
\int_{\R} F(u) \, \dx & \le \left(\frac{\eta}{4m} + \varepsilon\right) |u|_{2+4m}^{2+4m} + C |u|_{4+4m}^{4+4m}\\
& \le \left[\left(\frac{\eta}{4m} + \varepsilon\right) C_{2+4m}^{2+4m} \rho^{2m} + C C_{4+4m}^{4+4m} \rho^{1+2m-1/(2m)} |u^{(m)}|_2^{1/m}\right] |u^{(m)}|_2^2\\
& \le \left[\left(\frac{\eta}{4m} + \varepsilon\right) C_{2+4m}^{2+4m} \rho^{2m} + C C_{4+4m}^{4+4m} \rho^{1+2m-1/(2m)} \delta^{1/m}\right] |u^{(m)}|_2^2,
\end{align*}
which implies, in turn,
\begin{multline*}
8m J(u) = 4m |u^{(m)}|_2^2 - 8m \int_{\R} F(u) \, \dx\\
\ge \left[4m - \left(2 \eta + 8m \varepsilon\right) C_{2+4m}^{2+4m} \rho^{2m} - 8m C C_{4+4m}^{4+4m} \rho^{1+2m-1/(2m)} \delta^{1/m}\right] |u^{(m)}|_2^2\\
= \left[2 \left(2m - \eta C_{2+4m}^{2+4m} \rho^{2m}\right) - 8m C_{2+4m}^{2+4m} \rho^{2m} \varepsilon - 8m C C_{4+4m}^{4+4m} \rho^{1+2m-1/(2m)} \delta^{1/m}\right] |u^{(m)}|_2^2.
\end{multline*}
The claim then holds true taking
\[
\varepsilon = \frac{2m - \eta C_{2+4m}^{2+4m} \rho^{2m}}{16m C_{2+4m}^{2+4m} \rho^{2m}} \quad \text{and} \quad \delta = \min\left\{1,\left(\frac{2m - \eta C_{2+4m}^{2+4m} \rho^{2m}}{16m C C_{4+4m}^{4+4m} \rho^{1+2m-1/(2m)}}\right)^m\right\}.
\]
Now let $u \in \cD \cap \cM$ and set $s = \delta / |u^{(m)}|_2$ and $v = s \star u$ so that $v \in \cD \cap H^m(\R)$ and $|v'|_2 = \delta$. Then, Lemma \ref{l:fiber} yields
\[
J(u) \ge J(v) \ge \frac{2m - \eta C_{2+4m}^{2+4m} \rho^{2m}}{8m} \delta^2. \qedhere
\]
\end{proof}

\begin{Lem}\label{l:coercive}
If (\ref{F0})--(\ref{F4}) and \eqref{eq:rho} hold, then $J$ is coercive over $\cD \cap \cM$.
\end{Lem}
\begin{proof}
We refer to the proofs of \cite[Lemma 2.4]{BiegMed} or \cite[Lemma 2.5 (iv)]{JL2020CVPD}, which are similar.
\end{proof}

We need the following version of Lions's Lemma (cf. \cite[Lemma 3.1]{NonradMed}):
\begin{Lem}\label{l:Lions}
Let $G \in \cC(\R)$ such that $G(s) =o(s^2)$ as $u \to 0$. If $(u_n) \subset H^m(\R)$ is bounded and satisfies
\[
\lim_n \max_{y \in \R} \int_{y-r}^{y+r} u_n^2\, \dx=0
\]
for some $r>0$, then $\lim_{n \to +\infty} \int_{\R} |G(u_n)| \, \dx = 0$.
\end{Lem}

\begin{proof} Let $M>0$ be such that $\sup_n |u_n|_\infty < M$, and take any $\varepsilon > 0$ and $p>2$. Then, we find
$0 < \delta < M$ and $c_\varepsilon > 0$ such that
\begin{equation}
\begin{cases}
|G(s)| \leq \varepsilon s^2 &\text{ if } |s| \in [0, \delta],\\
|G(s)| \leq c_\varepsilon |s|^p &\text{ if } |s| \in  (\delta,M].
\end{cases}
\end{equation}
Hence, in view of Lions' lemma \cite[Lemma I.1]{Lions2}
, we get
$$
\limsup_{n\to\infty} \int_{\R^2} |G(u_n)| \, \dx \leq \varepsilon \limsup_{n\to\infty} \int_{\R^2} u_n^2 \, \dx.
$$
Letting $\varepsilon\to 0^+$, we conclude.
\end{proof}
Using Lemma \ref{l:Lions} and proceeding as in \cite[Proof of Theorem 1.4]{NonradMed}, we have the following result in the spirit of \cite{Gerard}. 
\begin{Prop}\label{p:split}
Suppose that $(u_n)\subset H^m(\R)$ is bounded. Then, there are sequences $(\widetilde{u}_i)_{i=0}^\infty \subset H^m(\R)$, $(y_n^i)_{i=0}^\infty \subset \R$ for every $n$, such that $y_n^0=0$, $|y_n^i-y_n^j|\to\infty$ as $n\to\infty$ for $i\neq j$, and passing to a subsequence, the following conditions hold for every $i\geq 0$:
$$
u_n(\cdot+y_n^i)\rightharpoonup \widetilde{u}_i \text{ in } H^m(\R) \text{ as } n\to\infty,
$$
\begin{equation}\label{eq11}
\lim_{n\to\infty} \int_\R |u_n^{(m)}|^2  \, \dx = \sum_{j=0}^i \int_\R |\widetilde{u}_j^{(m)}|^2 \, \dx + \lim_{n\to\infty} \int_\R |(v^i_n)^{(m)}|^2 \, \dx,
\end{equation}
where $v_n^i := u_n - \sum_{j=0}^i \widetilde{u}^j(\cdot-y_n^j)$, and
\begin{equation}\label{eq22}
\limsup_{n\to\infty} \int_\R G(u_n)\, \dx = \sum_{j=0}^\infty \int_\R G(\widetilde{u}_j) \, \dx
\end{equation}
for any function $G\colon \R\to[0,\infty)$ of class $\cC^1$ such that $|G'(s)|=\cO(|s|)$ as $s\to 0$ and $G(s)=o(s^2)$ as $s\to 0$.
\end{Prop}

\begin{Lem}\label{l:minimum}
If (\ref{F0})--(\ref{F4}) and \eqref{eq:rho} hold, then $\inf_{\cD \cap \cM} J$ is attained.
\end{Lem}
\begin{proof}
Let $(u_n)_n \subset \cD \cap \cM$ such that $\lim_n J(u_n) = \inf_{\cD \cap \cM} J$. From Lemma \ref{l:coercive}, we know that $u_n$ is bounded in $H^m(\R)$. Moreover, Proposition \ref{p:split} applied to the function $H$ implies the existence of the sequences $(\widetilde{u}_i)_{i=0}^{+\infty} \subset H^m(\R)$ satisfying \eqref{eq11} and \eqref{eq22}. Our first claim is that
\begin{equation}\label{claim11}
\exists  i \in \mathbb{N} \text{ such that }  0<\int_\R| \widetilde{u}_i^{(m)}|^2 \, \dx \leq \frac{1}{2m}\int_\R H(\widetilde{u}_i) \, \dx.
\end{equation}
Let $I := \{i \in \mathbb{N} : \widetilde{u}_i \ne 0\}$. In view of Lemma \ref{l:grad}, \eqref{eq22}, and the fact that $u_n \in \cM$, it is clear that $I \ne \emptyset$.
Arguing by contradiction, assume that
$$\int_\R| \widetilde{u}_i^{(m)}|^2 \, \dx >\frac{1}{2m}\int_\R H(\widetilde{u}_i) \, \dx \quad \forall i\in I.$$
Then, \eqref{eq11} and \eqref{eq22} imply that
\begin{align*}
\limsup_{n\to +\infty} \frac{1}{2m}\int_\R H(u_n) \, \dx &=\limsup_{n\to +\infty} \int_\R |u_n^{(m)}|^2 \, \dx \geq \sum_{j=0}^{+\infty} \int_\R| (\widetilde{u}_j)^{(m)}|^2 \, \dx\\
& =\sum_{j\in I}\int_\R| (\widetilde{u}_j)^{(m)}|^2 \, \dx
>\sum_{j\in I}\frac{1}{2m}\int_\R H (\widetilde{u}_j) \, \dx\\
& =\limsup_{n\to +\infty} \frac{1}{2m}\int_\R H(u_n) \, \dx,
\end{align*}
which is impossible. Consequently, \eqref{claim11} holds for some $i\in I$, and note that $r_i \ge 1$, where
\[
r_i^2 = \frac{\int_{\R} H(\widetilde{u}_i) \, \dx}{2m |\widetilde{u}_i^{(m)}|_2^2}.
\]
Additionally, a straightforward computation shows that $u := \widetilde{u}_i(r_i\cdot) \in \cM$. In fact, $u \in \cD \cap \cM$ because $r_i \ge 1$, therefore, from (\ref{F4}) and Fatou's Lemma, we have
\begin{align*}
\inf_{\cD \cap \cM} J & \le J(u) = \int_{\R} \frac{1}{4m} H(u) - F(u) \, \dx = \frac{1}{r_i} \int_{\R} \frac{1}{4m} H(\widetilde{u}_i) - F(\widetilde{u}_i) \, \dx\\
& \le \int_{\R} \frac{1}{4m} H(\widetilde{u}_i) - F(\widetilde{u}_i) \, \dx \\&\leq \liminf_{n \to +\infty} \int_{\R} \frac{1}{4m} H\bigl(u_n(x + y_{n}^i)\bigr) - F\bigl(u_n(x + y_{n}^i)\bigr) \, \dx\\
& \le \lim_{n \to +\infty} \int_{\R} \frac14 H(u_n) - F(u_n) \, \dx= \lim_{n \to +\infty} J(u_n) = \inf_{\cD \cap \cM} J,
\end{align*}
which implies that $r_i = 1$ and $u = \widetilde{u}_i$ minimises $J$ over $\cD \cap \cM$.
\end{proof}

\begin{proof}[Proof of Theorem \ref{th:mainSUP}]
The first part follows from Lemma \ref{l:minimum}. Now, let $u$ be the minimiser of $J$ over $\cD \cap \cM$ given therein. Our first claim is that the functional $\bigl(\Phi'(v),M'(v)\bigr) \colon H^m(\R) \to \R^2$ is surjective for every $v \in \cS \cap \cM$, where $\Phi(v) := |v|_2^2$ and
\[
M(v) := \int_{\R} \Bigl(|v^{(m)}|^2 - \frac{1}{2m} H(v)\Bigr) \, \dx.
\]
Indeed, given $v \in \cS \cap \cM$, we consider the curve $(0,+\infty)\ni s \mapsto s\star v \in\cS$, and the function $\psi_v(s)= M(s\star v)=s \varphi'_v(s)$. We notice that the curve $ s \mapsto s\star v$ is not tangent at $v$ to the manifold $\cM$, since otherwise we would obtain $\psi'_v(1)=\varphi''_v(1)=0$, in contradiction with Remark \ref{r:empty}). Thus, the manifolds $\cS$ and $\cM$ do not have the same tangent plane at $v$. Hence, from \cite[Proposition A.1]{MedSc}, there exist $\lambda \ge 0$ and $\theta \in \R$ such that
\begin{equation}\label{eq:LaThe}
(1 + \theta) \Bigl(-\frac{\mathrm{d}^2}{\dx^2}\Bigr)^{m} u + \lambda u = F'(u) + \frac{\theta}{4m} H'(u).
\end{equation}
If $\theta = -1$, then (\ref{F3}) and (\ref{F4}), together with \eqref{eq:LaThe}, imply
\begin{align*}
0 \le \lambda \int_{\R} |u|^2 \, \dx & = \int_{\R} \biggl(F'(u) u - \frac{1}{4m} H'(u) u\biggr) \, \dx < \int_{\R} \biggl(F'(u) u - \Bigl(\frac{1}{2m} + 1\Bigr) H(u) \biggr) \, \dx\\
& = \int_{\R} \biggl(2 F(u) - \frac{1}{2m} H(u)\biggr) \, \dx \le 0,
\end{align*}
a contradiction. As a consequence, from Proposition \ref{pr:Po}, $u$ satisfies also the Poho\v{z}aev and Nehari identity associated with \eqref{eq:LaThe}, whence
\[
(1 + \theta) \int_{\R} |u^{(m)}|^2 \, \dx = \frac{1}{2m} \int_{\R} H(u) + \frac{\theta}{4m} \bigl(H'(u) u - 2H(u)\bigr) \, \dx.
\]
Since $u \in \cM$, we obtain
\[
\theta \int_{\R} \bigl(H'(u) u - (2+4m) H(u)\bigr) \, \dx = 0,
\]
whence $\theta = 0$ in view of (\ref{F3}) and \eqref{eq:strict}. This proves that $\bigl(-\frac{\mathrm{d}^2}{\dx^2}\bigr)^{m} u + \lambda u = F'(u)$. Since $\lambda = 0$ if $u \in \cD \setminus \cS$, we only need to verify that $\lambda > 0$. Indeed, if by contradiction, $\lambda = 0$, then the Poho\v{z}aev identity yields
\[
\int_{\R} F(u) \, \dx = -\frac12 \int_{\R} |u^{(m)}|^2 \, \dx  < 0,
\]
which contradicts (\ref{F4}).
\end{proof}

\section*{Acknowledgements}
This work was supported by the Thematic Research Programme ``Variational and geometrical methods in partial differential equations'', University of Warsaw, Excellence Initiative - Research University.
J.S. is a member of GNAMPA (INdAM) and was supported by the programme \textit{Borse di studio per l'estero dell'Istituto Nazionale di Alta Matematica} and the GNAMPA project {\em Metodi variazionali e topologici per alcune equazione di Schr\"odinger nonlineari}.
The research project of Panayotis Smyrnelis is implemented in the framework of H.F.R.I call 
``Basic research Financing (Horizontal support of all Sciences)'' under the National Recovery and Resilience Plan ``Greece 2.0'' funded by the European Union - NextGenerationEU (H.F.R.I. Project Number: 016097).

\end{document}